\documentclass[a4paper, 10pt]{article}

\usepackage[T1]{fontenc}
\usepackage[latin1]{inputenc}
\usepackage{lmodern}  
\usepackage[english]{babel}
\usepackage{csquotes} 

\usepackage{textcomp} 
\usepackage{graphicx} 
\usepackage{amsmath} 
\usepackage{amsthm} 
\usepackage{amssymb} 
\usepackage{amsfonts}
\usepackage{mathtools}
\usepackage{thmtools} 
\usepackage{mathrsfs} 
\usepackage{mathdots} 
\usepackage{wasysym} 
\usepackage{extpfeil} 
\usepackage{verbatim} 
\usepackage{longtable} 
\usepackage{bigstrut} 
\usepackage[table, dvipsnames]{xcolor} 
\usepackage{floatpag} 
\usepackage{enumitem} 

\usepackage{xparse} 
\usepackage{mathscinet} 
\usepackage{adjustbox} 
\usepackage[obeyDraft]{todonotes} 

\usepackage[final,                  
            colorlinks,             
            plainpages = false,     
            linktoc=all,            
            linkcolor=black,        
            citecolor=black,         
            urlcolor=black,          
            filecolor=black         
            ]{hyperref}

\usepackage[citestyle=alphabetic,   
            bibstyle=alphabetic,    
            maxbibnames=50,         
            maxcitenames=3,         
            autocite=inline,        
            block=space,            
            hyperref=true,          
            backref=false,           
            backrefstyle=three,     
            date=short,             
            arxiv=pdf,              
            isbn=false,             
            url=false,              
            doi=false,               
            eprint=true,            
            giveninits=true,        
            sorting=nyt,            
            backend=bibtex8,
            ]{biblatex}

\usepackage{url}

\bibliography{cellular}

\hypersetup{pdftitle = {Cellularity of the p-Canonical Basis for Symmetric Groups},
            pdfauthor = {Lars Thorge Jensen},
            pdfkeywords = {p-canonical basis} {Hecke algebra} {cellular basis},
            pdfdisplaydoctitle      
            }
            
\usepackage[arrow, matrix, curve]{xy} 
\usepackage{tikz} 
\usetikzlibrary{arrows, arrows.meta, bending, fit, intersections, calc, external, shapes.misc, shapes.geometric, decorations.markings, decorations.pathreplacing}
\usepackage{pgfplots}
\pgfplotsset{compat=1.14}

\usepackage{genyoungtabtikz} 

\allowdisplaybreaks 

\setkeys{Gin}{draft=false}
            
\declaretheoremstyle[
spaceabove=\topsep, spacebelow=\topsep,
headfont=\normalfont\bfseries,
notefont=\mdseries, notebraces={(}{)},
bodyfont=\itshape,
postheadspace=\newline
]{break}

\declaretheoremstyle[
spaceabove=\topsep, spacebelow=\topsep,
headfont=\normalfont\bfseries,
notefont=\mdseries, notebraces={}{},
bodyfont=\itshape,
postheadspace=\newline
]{refbreak}

\declaretheorem[title=Theorem, style=plain, numberwithin=section]{thm}
\declaretheorem[title=Proposition, style=plain, numberlike=thm]{prop}
\declaretheorem[title=Lemma, style=plain, numberlike=thm]{lem}
\declaretheorem[title=Corollary, style=plain, numberlike=thm]{cor}

\declaretheorem[title=Assumption, style=plain, numberlike=thm]{ass}

\declaretheorem[title=Theorem, style=break, numberlike=thm]{thmlab}

\declaretheorem[title=Conjecture, style=break, numberlike=thm]{conj}

\declaretheorem[title=Definition, style=definition, numberlike=thm]{defn}

\declaretheorem[title=Example, style=remark, numberlike=thm]{ex}

\usepackage[capitalize]{cleveref} 
                          
\crefname{thm}{Theorem}{Theorems}
\crefname{prop}{Proposition}{Propositions}
\crefname{lem}{Lemma}{Lemmata}
\crefname{cor}{Corollary}{Corollaries}
\crefname{rem}{Reminder}{Reminders}
\crefname{defn}{Definition}{Definitions}

\crefname{thmlab}{Theorem}{Theorems}
\crefname{proplab}{Proposition}{Propositions}
\crefname{lemlab}{Lemma}{Lemmata}
\crefname{corlab}{Corollary}{Corollaries}
\crefname{remlab}{Reminder}{Reminders}
\crefname{conj}{Conjecture}{Conjectures}

\crefname{thmreflab}{Theorem}{Theorems}
\crefname{propreflab}{Proposition}{Propositions}
\crefname{lemreflab}{Lemma}{Lemmata}
\crefname{correflab}{Corollary}{Corollaries}
\crefname{remreflab}{Reminder}{Reminders}
\crefname{conjref}{Conjecture}{Conjectures}

\crefname{remark}{Remark}{Remarks}
\crefname{claim}{Claim}{Claims}
\crefname{ex}{Example}{Examples}

\crefname{section}{Section}{Sections}
\crefname{figure}{Figure}{Figures}
\crefname{equation}{}{}
\crefname{ass}{Assumption}{Assumptions}

\CompileMatrices

\entrymodifiers={+!!<0pt,\fontdimen22\textfont2>} 

\def\clap#1{\hbox to 0pt{\hss#1\hss}}

\makeatletter
\def\underbracket{%
    \@ifnextchar[{\@underbracket}{\@underbracket [\@bracketheight]}%
}
\def\@underbracket[#1]{%
    \@ifnextchar[{\@under@bracket[#1]}{\@under@bracket[#1][0.4em]}%
}
\def\@under@bracket[#1][#2]#3{
    \mathop{\vtop{\m@th \ialign {##\crcr $\hfil \displaystyle {#3}\hfil $%
    \crcr \noalign {\kern 3\p@ \nointerlineskip }\upbracketfill {#1}{#2}
    \crcr \noalign {\kern 3\p@ }}}}\limits}

\def\upbracketfill#1#2{$\m@th \setbox \z@ \hbox {$\braceld$}
    \edef\@bracketheight{\the\ht\z@}\bracketend{#1}{#2}
    \leaders \vrule \@height #1 \@depth \z@ \hfill
    \leaders \vrule \@height #1 \@depth \z@ \hfill \bracketend{#1}{#2}$}

\def\bracketend#1#2{\vrule height #2 width #1\relax}
\makeatother

\makeatletter
\def\thmt@refnamewithcomma #1#2#3,#4,#5\@nil{%
  \@xa\def\csname\thmt@envname #1utorefname\endcsname{#3}%
  \ifcsname #2refname\endcsname
    \csname #2refname\expandafter\endcsname\expandafter{\thmt@envname}{#3}{#4}%
  \fi
}
\makeatother

\makeatletter
\newcommand*\rel@kern[1]{\kern#1\dimexpr\macc@kerna}
\newcommand*\widebar[1]{%
  \begingroup
  \def\mathaccent##1##2{%
    \rel@kern{0.8}%
    \overline{\rel@kern{-0.8}\macc@nucleus\rel@kern{0.2}}%
    \rel@kern{-0.2}%
  }%
  \macc@depth\@ne
  \let\math@bgroup\@empty \let\math@egroup\macc@set@skewchar
  \mathsurround\z@ \frozen@everymath{\mathgroup\macc@group\relax}%
  \macc@set@skewchar\relax
  \let\mathaccentV\macc@nested@a
  \macc@nested@a\relax111{#1}%
  \endgroup
}
\makeatother

\makeatletter
\newcommand{\subjclass}[2][1991]{%
  \let\@oldtitle\@title%
  \gdef\@title{\@oldtitle\footnotetext{#1 \emph{Mathematics subject classification.} #2}}%
}
\newcommand{\keywords}[1]{%
  \let\@@oldtitle\@title%
  \gdef\@title{\@@oldtitle\footnotetext{\emph{Key words and phrases.} #1.}}%
}
\makeatother

\makeatletter
\newcommand{\extp}{\@ifnextchar^\@extp{\@extp^{\,}}}
   \def\@extp^#1{\mathop{\bigwedge\nolimits^{\!#1}}}
\makeatother

\DeclareMathOperator{\Mat}{Mat}
\DeclareMathOperator{\Hom}{Hom}
\DeclareMathOperator{\End}{End}

\DeclareMathOperator{\tr}{T}

\newcommand{\id}{\mathrm{id}}
\DeclareMathOperator{\sgn}{sgn} 
\DeclareMathOperator{\triv}{triv} 
\DeclareMathOperator{\geom}{geom} 
\DeclareMathOperator{\ch}{ch}

\DeclareMathOperator{\f}{f}

\DeclareMathOperator{\BS}{BS}
\DeclareMathOperator{\Spec}{Spec}
\DeclareMathOperator{\Irr}{Irr}

\newcommand{\defeq}{\ensuremath{\coloneqq}}

\newcommand{\cat}[1]{\ensuremath{\mathcal{#1}}}

\newcommand{\karalg}[1]{\ensuremath{\cat{K}ar(#1)}}

\newcommand{\G}[1]{\ensuremath{[\cat{#1}]}}
\newcommand{\Galg}[1]{\ensuremath{[#1]}}

\newcommand{\Z}{\ensuremath{\mathbb{Z}}}
\newcommand{\N}{\ensuremath{\mathbb{N}}}
\newcommand{\Q}{\ensuremath{\mathbb{Q}}}
\newcommand{\R}{\ensuremath{\mathbb{R}}}
\newcommand{\C}{\ensuremath{\mathbb{C}}}
\newcommand{\F}[1]{\ensuremath{\mathbb{F}_{#1}}}
\newcommand{\K}{\ensuremath{\mathbb{K}}}

\newcommand{\p}[1]{\ensuremath{\prescript{p}{}{#1}}}
\newcommand{\pre}[2]{\ensuremath{{}^{#1} #2}}

\newcommand{\heck}{\ensuremath{\mathcal{H}}}
\newcommand{\heckMod}[1]{\ensuremath{\mathcal{H}(#1)}}
\newcommand{\weylMod}[2][\C]{\ensuremath{#1\W(#2)}}

\newcommand{\std}[2][H]{\ensuremath{#1_{#2}}}
\newcommand{\kl}[2][H]{\ensuremath{\underline{#1}_{#2}}}
\DeclareDocumentCommand{\pkl}{O{p} O{H} m}{\ensuremath{\prescript{#1}{}{\underline{#2}}_{#3}}}
\newcommand{\pfpelt}[2][p]{{}^{#1}b_{\mathbf{#2}}}
\newcommand{\desc}[1]{\ensuremath{\mathcal{#1}}}

\newcommand{\expr}[1]{\ensuremath{\underline{#1}}}

\newcommand{\cle}[2][p]{\ensuremath{\, \overset{#1}{\underset{#2}{\leqslant}} \,}}
\newcommand{\cge}[2][p]{\ensuremath{\, \overset{#1}{\underset{#2}{\geqslant}} \,}}

\newcommand{\clt}[2][p]{\ensuremath{\, \overset{#1}{\underset{#2}{<}} \,}}

\newcommand{\ceq}[2][p]{\ensuremath{\, \overset{#1}{\underset{#2}{\sim}} \,}}

\newcommand{\T}{\ensuremath{\cat{T}}}

\newcommand{\rt}[1][]{\ensuremath{\alpha_{#1}}}
\newcommand{\cort}[1][]{\ensuremath{\alpha_{#1}^{\vee}}}

\newcommand{\charlat}{\ensuremath{X}}
\newcommand{\cocharlat}{\ensuremath{X^{\vee}}}

\newcommand{\Wid}{\ensuremath{e}}

\newcommand{\W}{\ensuremath{W}}

\newcommand{\pafp}[2][p]{{}^{#1} \!\mkern 1mu {\mathbf a}_{\mathbf{#2}}} 

\newcommand{\HC}[1][]{\ensuremath{\mathbf{H}^{#1}}}

\newcommand{\Homnl}[2][]{
   \ifthenelse{ \equal{#1}{} } { \ensuremath{\Hom_{\nless #2}} } 
   { \ensuremath{\Hom_{\nless #2, #1}} } }
\newcommand{\Endnl}[2][]{
   \ifthenelse{ \equal{#1}{} } { \ensuremath{\End_{\nless #2}} } 
   { \ensuremath{\End_{\nless #2, #1}} } }
\newcommand{\D}{\ensuremath{\mathbb{D}}}

\tikzset{%
  DynNode/.style={circle, inner sep=2pt, draw=black, fill=white},
  Greater/.style={pos=0.65, inner sep=0mm, outer sep=0mm},
  highlight/.style={rectangle,rounded corners,fill=red!15,draw=red,
    fill opacity=0.5,thick},
  root/.style={draw, color=black, thick, ->},
  plane/.style={draw, color=black, very thin},
  origin/.style={fill, color=black},
  sline/.style={color=Red, thick},
  tline/.style={color=NavyBlue, thick},
  uline/.style={color=Goldenrod, thick},
  bendBelow/.style={bend left=70, looseness=2},
  bendAbove/.style={bend right=70, looseness=2},
  object/.style={circle, fill, inner sep=1.5pt, outer sep=0mm},
  labelling/.style={outer sep=0mm, inner sep=0mm},
  1morph/.style={->, shorten >= 0.5pt, >=stealth'},
  2morph/.style={-implies,double,double equal sign distance,
                 shorten >=2pt, shorten <=3pt},
  spot/.style={color=black, thin, dashed},
  s/.style={color=Red},
  t/.style={color=NavyBlue},
  u/.style={color=Goldenrod},
  line/.style={draw, line width=2pt},
  dot/.style={fill, thin},
  sph/.style={fill, color=black!20, opacity=0.5},
  squig/.style={decoration=snake, decorate, ->},
  root/.style={draw, color=black, line width=2pt, ->},
  myptr/.style={decoration={markings,mark=at position 1 with %
    {\arrow[scale=3,>=stealth]{>}}},postaction={decorate}},
  on each segment/.style={
    decorate,
    decoration={
      show path construction,
      moveto code={},
      lineto code={
        \path [#1]
        (\tikzinputsegmentfirst) -- (\tikzinputsegmentlast);
      },
      curveto code={
        \path [#1] (\tikzinputsegmentfirst)
        .. controls
        (\tikzinputsegmentsupporta) and (\tikzinputsegmentsupportb)
        ..
        (\tikzinputsegmentlast);
      },
      closepath code={
        \path [#1]
        (\tikzinputsegmentfirst) -- (\tikzinputsegmentlast);
      },
    },
  },
  mid arrow/.style={postaction={decorate,decoration={
        markings,
        mark=at position .5 with {\arrow[#1]{stealth}}
      }}},
  sline/.style={draw, line width=1pt, postaction={on each segment={mid arrow=black}}},
  sregion/.style={fill, opacity=0.2},
  st/.style={fill=Fuchsia},
  su/.style={fill=YellowOrange},
  tu/.style={fill=ForestGreen},
  clabel/.style={fill=none, red}, 
  str/.style={<->}
}

\begin{document}

\title{Cellularity of the $p$-Canonical Basis for Symmetric Groups}
\author{Lars Thorge Jensen}
\date{} 

\maketitle

\begin{abstract}
   For symmetric groups we show that the $p$-canonical basis can be extended 
   to a cell datum for the Iwahori-Hecke algebra $\heck$ and that the two-sided 
   $p$-cell preorder coincides with the Kazhdan-Lusztig two-sided cell preorder.
   Moreover, we show that left (or right) $p$-cells inside the same two-sided 
   $p$-cell for Hecke algebras of finite crystallographic Coxeter systems are 
   incomparable (Property A).
\end{abstract}

\tableofcontents

\section{Introduction}

The Hecke algebra of a crystallographic Coxeter system admits several 
geometric or algebraic categorifications (see \cite{KL, EW2}). 
The canonical bases arising from these categorifications coincide with the 
famous Kazhdan-Lusztig basis (see \cite{KL2, EW1}) in the characteristic
$0$ setting and give rise to the $p$-canonical or $p$-Kazhdan-Lusztig basis
in the characteristic $p > 0$ setting (see \cite{JW, RWTiltPCan}). 


The study of cells with respect to the $p$-canonical basis of the Hecke algebra 
was initiated in \cite{JeABC} and continued in \cite{JPInd}. In this paper,
we apply the Perron-Frobenius theorem to $p$-cells and to tackle several questions 
about the various $p$-cell preorders. One of the most important features of the 
Kazhdan-Lusztig cell preorders is Property A which states that left (or right) 
Kazhdan-Lusztig cells in the same two-sided Kazhdan-Lusztig cell are incomparable 
with respect to the left (or right) cell preorder. We will prove Property A of 
$p$-cells for finite crystallographic Coxeter groups in this paper following ideas 
of \cite{KM}.

Along the way, we introduce $p$-special modules and $p$-families which 
are generalizations of Lusztig's special modules and families of irreducible representations of a Weyl group to the $p$-canonical basis. These concepts were originally
introduced by Lusztig in \cite{LuUnipRepsOfE8, LuSpecialReps, LuSpecialRepsII} 
and have played an important role in determining the complex irreducible characters
of finite reductive groups (via character sheaves). The connections of $p$-special
modules and $p$-families to the representation theory of finite reductive groups
are completely unclear to the author and thus merit further study.

In the last part of the paper, we study consequences of Property A for $p$-cells
of symmetric groups. One of the main results of \cite{JeABC} is the characterization 
of $p$-cells for symmetric groups in terms of the Robinson-Schensted correspondence 
which gives a bijection $w \mapsto (P(w), Q(w))$ between the symmetric group $S_n$ 
and pairs of standard tableaux of the same shape with $n$ boxes. As a consequence 
$p$-cells and Kazhdan-Lusztig cells coincide for all primes $p$: The two-sided cell 
of $w \in S_n$ is given by the set of elements in $S_n$ whose $Q$-symbols have the 
same shape as $Q(w)$. In other words, two-sided cells of $S_n$ are in bijection 
with partitions of $n$.

Even though the cells coincide, it was an open question to relate the two-sided 
$p$-cell preorder and the Kazhdan-Lusztig two-sided cell preorder on the level of 
two-sided cells. As a special case of the Lusztig-Vogan bijection (which was proven
in \cite{BLuBij}) the Robinson-Schensted correspondence gives 
an order-preserving bijection between the set of two-sided cells equipped with the 
Kazhdan-Lusztig two-sided cell preorder and the set of partitions equipped with 
the dominance order. We extend this result to the two-sided $p$-cell preorder 
showing that the two-sided $p$-cell preorder and the Kazhdan-Lusztig two-sided 
cell preorder coincide for all primes $p$.

When Graham and Lehrer introduced cellular algebras in \cite{GLCellAlgs}, the
Hecke algebra of a symmetric group together with the Kazhdan-Lusztig basis was 
one of the motivating examples. Finally, we show that the $p$-canonical basis
can also be extended to a cell datum for the Hecke algebra of a symmetric group.
Let's recall the main result of \cite[Theorem 1.1]{LMNRightCell}:
\begin{thm}
   Let $x, y \in S_n$ with $x \leqslant y$. Then there exist $N \geqslant n$
   and $v, w \in S_N$ with $v \leqslant w$ such that
   \begin{enumerate}
      \item $v$ and $w$ belong to the same right cell,
      \item the singularity of $X_w$ at $v$ is smoothly equivalent to the
            singularity of $X_y$ at $x$.
   \end{enumerate}
\end{thm}

As explained in \cite[\S4.3]{LMNRightCell} this result allows to embed
any example from \cite{WTorsion} in which the Kazhdan-Lusztig and the $p$-canonical
basis differ into the same cell in some larger symmetric group and thus
produces many examples where the $p$-canonical basis gives interesting
basis of Specht modules. This adds to the interest in this newly defined cell
datum for the Hecke algebra of a symmetric group.

\subsection{Structure of the Paper}

\begin{description}
   \item[\Cref{secBack}] We introduce notation for crystallographic Coxeter systems
         and their Hecke algebras. Then we recall important results about
         the diagrammatic category of Soergel bimodules and the $p$-canonical basis.
         Finally, we remind the reader of the Perron-Frobenius Theorem.
   \item[\Cref{secSpecialMods}] We apply the Perron-Frobenius Theorem to $p$-cells 
         to introduce $p$-special modules and $p$-families and to prove some of their
         elementary properties.
   \item[\Cref{secPropertyA}] We prove Property A for finite Weyl groups, 
         showing that left $p$-cells within the same two-sided $p$-cell
         are incomparable with respect to the left $p$-cell preorder.
   \item[\Cref{secPKLCellular}] This section deals with some consequences for $p$-cells
         of symmetric groups. First, we recall the characterization of $p$-cells in
         terms of the Robinson-Schensted correspondence. Then we show that the
         two-sided $p$-cell preorder corresponds to the dominance order on partitions.
         Finally, we prove that the $p$-canonical basis can be extended to a cellular 
         datum for the Hecke algebra.
\end{description}

\subsection{Acknowledgements}
Since this project goes back to a question by Peter McNamara at the conference
on ``New Connections in Representation Theory'' in Mooloolaba, Australia, in 
February 2020, I would like to thank him for asking me the question and later
encouraging me to write up the answer. Moreover, I am grateful to the 
Henri Poincar\'e Institute in Paris for the good working conditions and 
to Simon Riche and the ERC project ModRed for the financial support 
for my research visit in Sydney. I would also like to thank Leonardo Patimo 
and Shotaro Makisumi for helpful discussions. The author has received funding
from the European Research Council (ERC) under the European Union's Horizon 2020 
research and innovation programme (grant agreement No 677147).

\section{Background}
\label{secBack}

\subsection{Crystallographic Coxeter Systems}

Let $S$ be a finite set and $(m_{s, t})_{s, t \in S}$ be a matrix with entries in
$\N \cup \{\infty\}$ such that $m_{s,s} = 1$ and $ m_{s, t} = m_{t, s} \geqslant 2$ 
for all $s \neq t \in S$. Denote by $W$ the group generated by $S$ subject to the relations 
$(st)^{m_{s, t}} = 1$ for $s, t \in S$ with $m_{s, t} < \infty$. We say that $(W, S)$
is a \emph{Coxeter system} and $W$ is a \emph{Coxeter group}. The Coxeter group $W$
comes equipped with a length function $l: W \rightarrow \N$ and the Bruhat order
$\leqslant$ (see \cite{Hum} for more details). A Coxeter system $(W, S)$ is 
called \emph{crystallographic} if $m_{s, t} \in \{2, 3, 4, 6, \infty\}$ for 
all $s \neq t \in S$. We denote the identity of $W$ by $\Wid$. For $w \in W$ we
define its \emph{left descent set} via
\[ \desc{L}(w) \defeq \{ s \in S \; \vert \; l(sw) < l(w) \} \text{.} \]
The \emph{right descent set} of $w$ is given by $\desc{R}(w) \defeq \desc{L}(w^{-1})$.

From now on, fix a generalized Cartan matrix $A = (a_{i,j})_{i, j \in J}$ (see 
\cite[\S1.1]{TiGrpKM}). Let $(J, \charlat, \{ \rt[i] : i \in J\}, 
\{ \cort[i] : i \in J \})$ be an associated Kac-Moody root datum (see 
\cite[\S1.2]{TiGrpKM} for the definition). Then $X$ is a finitely generated 
free abelian group, and for $i \in J$ we have elements $\rt[i]$ and $\cort[i]$ 
of $X$ and $\cocharlat = \Hom_{\Z}(\charlat, \Z)$ respectively that satisfy
$a_{i,j} = \cort[i](\rt[j])$ for all $i, j \in J$.

To $A$ we associate a crystallographic Coxeter system $(W, S)$ as follows:
Choose a set of simple reflections $S$ of cardinality $\lvert J \rvert$ and
fix a bijection $S \overset{\sim}{\rightarrow} J$, $s \mapsto i_s$. For $s \ne t \in S$
we define $m_{s, t}$ to be $2$, $3$, $4$, $6$, or $\infty$ if $a_{i_s, i_t} a_{i_t, i_s}$
is $0$, $1$, $2$, $3$, or $\geqslant 4$ respectively. Most of the time, we will
work with a Cartan matrix, so that $W$ is a finite Weyl group.

Fix a commutative ring $k$. ${}^k V \defeq \cocharlat \otimes_{\Z} k$ 
yields a balanced, potentially non-faithful realization of the Coxeter system over $k$. 
Set ${}^k V^{\ast} \defeq \Hom_{k}({}^k V, k)$ and note that ${}^k V^{\ast}$ is 
isomorphic to $\charlat \otimes_{\Z} k$. A realization obtained in this way is called
a \emph{Cartan realization} (see \cite[\S10.1]{AMRWFreeMonodromicMixedTiltSheaves}).
Throughout, we will assume our realization to satisfy:
\begin{ass}[Demazure Surjectivity]
   \label{assDemazureSurj}
   The maps $\rt[s]: {}^k V \rightarrow k$ and $\cort[s]: {}^k V^{\ast} 
   \rightarrow k$ are surjective for all $s \in S$.
\end{ass}
\noindent This is automatically satisfied if $2$ is invertible in $k$ or if the
Coxeter system $(W, S)$ is of simply-laced type and of rank $\lvert S \rvert \geqslant 2$.

We denote by $R = S({}^k V^{\ast})$ the symmetric algebra of ${}^k V^{\ast}$
over $k$ and view it as a graded ring with ${}^k V^{\ast}$ in degree $2$. 
Given a graded $R$-bimodule $B = \bigoplus_{i\in \Z} B^i$, we denote by 
$B(1)$ the shifted bimodule with $B(1)^i = B^{i+1}$.

\subsection{The Hecke Algebra}

The Hecke algebra is the free $\Z[v, v^{-1}]$-algebra with $\{ \std{w} \; \vert \; 
w \in W \}$ as basis, called the \emph{standard basis}, and multiplication determined by:
\begin{alignat*}{2}
   \std{s}^2 &= (v^{-1} - v) \std{s} + 1  \qquad && \text{for all } s \in S \text{,} \\
   \std{x} \std{y} &= \std{xy} && \text{if } l(x) + l(y) = l(xy) \text{.}
\end{alignat*}

There is a unique $\Z$-linear involution $\widebar{(-)}$ on $\heck$ satisfying
$\widebar{v} = v^{-1}$ and $\widebar{\std{x}} = \std{x^{-1}}^{-1}$. The Kazhdan-Lusztig 
basis element $\kl{x}$ is the unique element in $\std{x} + \sum_{y < x} v\Z[v] H_y$ 
that is invariant under $\widebar{(-)}$. This is Soergel's normalization from 
\cite{S2} of a basis introduced originally in \cite{KL}.

Let $\iota$ be the $\Z[v, v^{-1}]$-linear anti-involution on $\heck$ satisfying 
$\iota(\std{s}) = \std{s}$ for $s \in S$ and thus $\iota(\std{x}) = \std{x^{-1}}$.

\subsection{The Diagrammatic Category of Soergel Bimodules}

In this section, we introduce the diagrammatic category of Soergel bimodules. 
The main reference for this is \cite{EW2} (see also \cite{E3} in the dihedral 
case and \cite{EKh} in type $A$).

Let $\BS$ be the diagrammatic category of Bott-Samelson bimodules
as introduced in \cite[\S2.3]{JW}. It is a diagrammatic, strict monoidal 
category enriched over $\Z$-graded left $R$-modules.

Let $\HC$ be the Karoubian envelope of the graded version of the additive closure
of $\BS$, in symbols $\HC = \karalg{\BS}$. We call $\HC$ the 
\emph{diagrammatic category of Soergel bimodules}. In other words, in the passage 
from $\BS$ to $\HC$ we first allow direct sums and grading shifts (restricting 
to degree preserving homomorphisms) and then the taking of direct summands. The 
following properties can be found in \cite[Lemma 6.24, Theorem 6.25 and 
Corollary 6.26]{EW2}:

\begin{thmlab}[Properties of $\HC$]
   \label{thmDiagProps}
   Let $k$ be a complete, local, integral domain (e.g. a field or the $p$-adic 
   integers $\Z_p$).
   \begin{enumerate}
      \item $\HC$ is a Krull-Schmidt category.
      \item For all $w \in W$ there exists a unique indecomposable object $\pre{k}{B}_w 
            \in \HC$ which is a direct summand of $\expr{w}$ for any reduced 
            expression $\expr{w}$ of $w$ and which is not isomorphic to a grading 
            shift of any direct summand of any expression $\expr{v}$ for $v < w$. 
            In particular, the object $\pre{k}{B}_w$ does not depend up to isomorphism 
            on the reduced expression $\expr{w}$ of $w$.
      \item The set $\{ \pre{k}{B}_w \; \vert \; w \in W\}$ gives a complete set of 
            representatives of the isomorphism classes of indecomposable objects 
            in $\HC$ up to grading shift.
      \item There exists a unique isomorphism of $\Z[v, v^{-1}]$-algebras
            \[ \ch: \Galg{\HC} \longrightarrow \heck \]
            sending $[\pre{k}{B}_s]$ to $\kl{s}$ for all $s \in S$, where $\Galg{\HC}$ 
            denotes the split Grothendieck group of $\HC$. \textup{(}We view $\Galg{\HC}$ as a
            $\Z[v, v^{-1}]$-algebra as follows: the monoidal structure
            on $\HC$ induces a unital, associative multiplication and $v$
            acts via $v[B] \defeq [B(1)]$ for an object $B$ of $\HC$.\textup{)}
   \end{enumerate}
\end{thmlab}

It should be noted that we do not have a diagrammatic presentation of $\HC$ 
as determining the idempotents in $\BS$ is usually extremely difficult. 

\subsection{The \texorpdfstring{$p$}{p}-canonical Basis and \texorpdfstring{$p$}{p}-Cells}

In this section, we recall the definition of the $p$-canonical basis
and its elementary properties (see \cite{JW, JeABC}). 
Let $k$  be a field of characteristic $p \geqslant 0$. Note that the $p$-canonical basis
depends on $p$, but not on the explicit choice of $k$. 

\begin{defn}
   Define $\pkl{w} = \ch([\pre{k}{B}_w])$ for all $w \in W$ where $\ch: \Galg{\HC} 
   \overset{\cong}{\longrightarrow} \heck$ is the isomorphism of $\Z[v, v^{-1}]$-algebras
   introduced earlier.
\end{defn}

We will frequently use the following elementary properties
of the $p$-canonical basis which can be found in \cite[Proposition 4.2]{JW}
unless stated otherwise:

\begin{prop}
   \label{propPCanProps}
   For all $x, y \in W$ we have:
   \begin{enumerate}
      \item $\widebar{\pkl{x}} = \pkl{x}$, i.e. $\pkl{x}$ is self-dual,
      \item $\iota(\pkl{x}) = \pkl{x^{-1}}$ and thus in particular 
            $\p{m}_{y, x} = \p{m}_{y^{-1}, x^{-1}}$ as well as
            $\pre{p}{h}_{y, x} = \pre{p}{h}_{y^{-1}, x^{-1}}$,
      \item $\pkl{x} \pkl{y} = \sum_{z \in W} \p{\mu}^{z}_{x, y} \pkl{z}$
            with self-dual $\p{\mu}^{z}_{x, y} \in \Z_{\geqslant 0}[v, v^{-1}]$,
      \item $\pkl{x} = \kl{x}$ for $p = 0$ (see \cite{EW1}) and $p \gg 0$ 
            (i.e. there are only finitely many primes for which $\pkl{x} \neq \kl{x}$).
   \end{enumerate}
\end{prop}

Let us recall the definition of $p$-cells (see \cite[\S3.1]{JeABC} for more details)
which is an obvious generalization of a notion introduced by Kazhdan-Lusztig in \cite{KL}:

\begin{defn}
   For $h \in \heck$ we say that \emph{$\pkl{w}$ appears with non-zero coefficient
   in $h$} if the coefficient of $\pkl{w}$ is non-zero when expressing $h$ in
   the $p$-canonical basis.
   
   Define a preorder $\cle{R}$ (resp. $\cle{L}$) on $W$ as follows:
   $x \cle{R} y$ (resp $x \cle{L} y$) if and only if $\pkl{x}$ appears with non-zero
   coefficient in $\pkl{y} h$ (resp. $h \pkl{y}$) for some $h \in \heck$.
   Define $\cle{LR}$ to be the preorder generated by $\cle{R}$ and $\cle{L}$, in 
   other words we have: $x \cle{LR} y$ if and only if $\pkl{x}$ appears with non-zero
   coefficient in $h \pkl{y} h'$ for some $h, h' \in \heck$.
   
   The left, right, or two-sided $p$-cells are the equivalence classes with respect 
   to the corresponding preorders respectively.
\end{defn}

Let $R$ be a fixed commutative ring with unit and let $A$ be an $R$-algebra with fixed
$R$-basis indexed by $B$. For a subset $J \subseteq B$ we denote by $A(J)$ the $R$-span 
of the basis elements indexed by elements in $J$. Let $C$ be a left $p$-cell in $\W$. 
To simplify notation, we write
\begin{align*}
   \cle{L} C &\defeq \{ x \in \W \; \vert \; x \cle{L} y \text{ for some } y \in C \} \\
   \clt{L} C &\defeq \{ x \in \W \; \vert \; x \cle{L} y \text{ for some } y \in C 
   \text{ and } y \notin C \}
\end{align*}
Unless stated otherwise, we will consider throughout the $p$-canonical basis
of $\heck$ (and its extension of scalars to $\C$). Recall that $C$ gives rise 
to a cell module
\[ \heckMod{C} = \heckMod{\cle{L} C} / \heckMod{\clt{L} C} \]
which is a left module for $\heck$ (similarly for right or two-sided $p$-cells).
The $p$-canonical basis elements indexed by elements in $C$ gives rise to a basis 
of $\heckMod{C}$.

\subsection{The Perron-Frobenius Theorem}

Our main technical tool will be the Perron-Frobenius theorem (see \cite{PeMat, Fr1, Fr2})
which was published over a $100$ years ago. A more modern exposition can be found in
\cite[Vol. 2, Chapter XIII]{GaMat} or in \cite[Kapitel IV]{Hupp}.

\begin{thm}[Perron-Frobenius]
   \label{thmPF}
   Let $M \in \Mat_{k\times k}(\R_{>0})$. Then there exists $\lambda \in \R_{>0}$,
   called the Perron-Frobenius eigenvalue of $M$, such that the following 
   statements holds:
   \begin{enumerate}
      \item $\lambda$ is an eigenvalue of $M$.
      \item Any other eigenvalue $\mu \in \C$ of $M$ satisfies $\lvert \mu \rvert < \lambda$,
            so $\lambda$ gives the spectral radius of $M$.
      \item The eigenvalue $\lambda$ has algebraic multiplicity $1$.
      \item There exists $v \in \R^k_{> 0}$ such that $Mv = \lambda v$. There exists also
            $\hat{v} \in \R^k_{> 0}$ such that $\hat{v}^{\tr} M = \lambda \hat{v}^{\tr}$.
      \item Any $w \in \R^k_{\geqslant 0}$ which is an eigenvector for $M$ is a scalar
            multiple of $v$ and similarly for $\hat{v}$.
      \item If $v$ and $\hat{v}$ are normalized such that $\hat{v}^{\tr} v = (1)$, then
            \[ \lim_{n \rightarrow \infty} \frac{M^n}{\lambda^n} = v\hat{v}^{\tr}\text{.} \]
   \end{enumerate}
\end{thm}

\section{\texorpdfstring{$p$}{p}-Families and \texorpdfstring{$p$}{p}-Special Modules}
\label{secSpecialMods}

In this section, we will apply some results of \cite{KM} to the 
$p$-canonical basis of the complex group ring of a finite Weyl group.
Thus we assume that $(W, S)$ is a finite Weyl group throughout the section.
For the sake of completeness, we give all the proofs.

Denote by $\C\W = \C \otimes_{\Z[v, v^{-1}]} \heck$ the scalar extension
of $\heck$ to $\C$ where we specialize $v$ to $1$. For the rest of the paper
fix a set of coefficients $\mathbf{c} = \{c_w\}_{w \in \W}$ with $c_w \in \R_{> 0}$.
$\mathbf{c}$ determines a Perron-Frobenius element
\[ \pfpelt{c} = \sum_{w \in \W} c_w \otimes \pkl{w} \in \C\W \text{.} \]

Let $C$ be a left $p$-cell in $\W$.  We will denote by 
$\weylMod{C} = \C \otimes_{\Z[v, v^{-1}]} \heckMod{C}$ the extension of scalars 
of the corresponding left cell module $\heckMod{C}$ to $\C$. As an immediate 
consequence of \cref{thmPF} we get in this setting (see \cite[Corollary 3]{KM}):
\begin{cor}
   The left (resp. right) action of $\pfpelt{c}$ on
   $\weylMod{C}$ gives a Perron-Frobenius eigenvalue $\pafp{c}(C)$.
   There is up to isomorphism a unique irreducible $\C\W$-module 
   $L_{C, \mathbf{c}}$ occuring as composition factor of the cell module $\weylMod{C}$
   such that $\pafp{c}(C)$ is afforded by the action of $\pfpelt{c}$ 
   on $L_{C, \mathbf{c}}$. Moreover, $[\weylMod{C} : L_{C, \mathbf{c}}]$ is equal to $1$.
\end{cor}

The following result shows that we can simplify our notation (see \cite[Theorem 5]{KM}):
\begin{thm}
   $L_{C, \mathbf{c}}$ does up to isomorphism not depend on the choice of $\mathbf{c}$.
\end{thm}
\begin{proof}
   For simplicity, denote by $n = \lvert W \rvert$ the cardinality of $W$.
   Consider the map $L_{C, -}: \R_{>0}^n \rightarrow \Irr(\C\W)$ which sends $\mathbf{d}$
   to $L_{C, \mathbf{d}}$. Equipping $\Irr(\C\W)$ with the discrete topology,
   we claim that this map is continuous. Since $R_{>0}^n$ is connected, 
   any continuous function to a discrete set has to be constant. Therefore,
   the claim implies the theorem.
   
   To prove the claim, we will show that the preimage $X_L \subseteq \R_{>0}^n$ of 
   any $L \in \Irr(\C\W)$ under this map is closed.
   Assume $X_L$ is non-empty and let $\mathbf{d}_i \in \R_{>0}^n \cap X_L$ be
   a sequence that converges to $\mathbf{d} \in \R_{>0}^n$. Let $L_1, L_2, \dots, 
   L_k = L$ be the list of simple subquotients of $\weylMod{C}$. 
   Denote by $M_{\mathbf{d}_i}$ the linear operator on $\weylMod{C}$ the
   element $\pfpelt{\mathbf{d}_{\mathrm{i}}} \in \C\W$ induces. As $d_i \in X_L$,
   we have for all $i \in \N$ by \cref{thmPF}:
   \begin{align*}
      \pafp{d_{\mathnormal{i}}}(C) &= 
      \max \{ \lvert \mu \rvert \; \vert \; \mu \in \Spec(M_{\mathbf{d}_i} \vert_{L}) \}\\
      &> \sup_{1 \leqslant j < k} \{ \lvert \mu \rvert \; \vert \; \mu \in \Spec(M_{\mathbf{d}_i} \vert_{L_j}) \}
   \end{align*}
   
   Therefore, we get in the limit as the spectrum of a matrix depends continuously
   on the matrix (see \cite[\S 5.2.3]{SeMats}):
   \begin{align*}
      \max \{ \lvert \mu \rvert \; \vert \; \mu \in \Spec(M_{\mathbf{d}} \vert_{L}) \}
      \geqslant \sup_{1 \leqslant j < k} \{ \lvert \mu \rvert \; \vert \; 
         \mu \in \Spec(M_{\mathbf{d}_i} \vert_{L_j}) \}
   \end{align*}
   Since the Perron-Frobenius eigenvalue has multiplicity one, it follows that 
   $L_{C, \mathbf{d}}$ is still isomorphic to $L$. It follows that $X_L$ is closed
   and thus finishes the proof of the claim (and the theorem).
\end{proof}

Due to the last result, we can and will drop $\mathbf{c}$ in the subscript and 
denote $L_{C, \mathbf{c}}$ by $L_C$ from now on. The following result shows that 
$L_C$ is an invariant of two-sided $p$-cells (see \cite[Theorem 6]{KM}):
\begin{thm}
   \label{thmInv}
   For any other left $p$-cell $C'$ which belongs to the same two-sided 
   $p$-cell as $C$, we have $L_C \cong L_{C'}$ and $\pafp{c}(C) = \pafp{c}(C')$.
\end{thm}
\begin{proof}
   Denote by $J$ the two-sided cell that contains $C$ and $C'$.
   We may assume that $C'$ is maximal with respect to $\cle{L}$ in $J$.

   First, we will construct a non-zero homomorphism 
   $\varphi: \weylMod{C} \longrightarrow \weylMod{C'}$ of $\C \W$-modules.
   For any $u \in C$ and $v \in C'$ there exist $h', h \in \heck$ such that
   $\pkl{v}$ occurs with non-zero coefficient in $h' \pkl{u} h$ by the definition
   of two-sided $p$-cells. It follows that $\pkl{u} h$ intersects the set 
   $\{ x \in \W \; | \; x \cge{L} v\} \cap J$ non-trivially. By our assumption of
   maximality the intersection $\{ x \in \W \; | \; x \cge{L} v\} \cap J$ 
   coincides with $C'$. Thus right multiplication by $h$ and projection onto
   $\weylMod{C'}$ defines the non-zero homomorphism $\varphi$. Observe 
   that $\varphi$ sends any linear combination of the $p$-canonical basis 
   of $\weylMod{C}$ with strictly positive coefficients to a non-zero
   linear combination of the $p$-canonical basis of $\weylMod{C'}$ with
   non-negative coefficients.
   
   By \cref{thmPF} (iv) there exists an eigenvector $v$ of $\pfpelt{c}$ on
   $\weylMod{C}$ with eigenvalue $\pafp{c}(C)$ that is a linear combination 
   of the p-canonical basis with strictly positive coefficients. It 
   follows that $\varphi(v)$ is non-zero and an eigenvector of $\pfpelt{c}$ 
   in $\weylMod{C'}$ with eigenvalue $\pafp{c}(C)$. Moreover,
   $\varphi(v)$ is a linear combination of the $p$-canonical basis of $\weylMod{C'}$ 
   with non-negative coefficients. Therefore, the corresponding eigenvalue
   is the Perron-Frobenius eigenvalue of $\pfpelt{c}$ on $\weylMod{C'}$
   by \cref{thmPF} (v). This implies $\pafp{c}(C) = \pafp{c}(C')$. As $v \in L_C$, 
   the simple subquotient $L_C$ is not annihilated by $\varphi$. Schur's 
   lemma implies that $\varphi$ induces an isomorphism $L_C \cong L_{C'}$.
\end{proof}

Kildetoft and Mazorchuk prove the following interesting results as well
(see \cite[Proposition 13]{KM}):
\begin{prop}
   \label{propFamilies}
   \begin{enumerate}
    \item The number of left $p$-cells in the two-sided $p$-cell
          of $C$ is given by $\dim L_C$.
    \item The two sided $p$-cells induce a partition of the irreducible
          representations of $\W$ as follows:
          \[ \Irr(\W) = \bigcup_{\substack{J \subseteq \W\\\text{two-sided }p\text{-cell}}} 
               \{ L \in \Irr(\W) \; \vert \; [ \weylMod{J}, L ] \neq 0 \} \]
          In particular, any simple subquotient of $\weylMod{C}$ different from $L_C$ is 
          not isomorphic to $L_{C'}$ for any other left $p$-cell $C'$.
   \end{enumerate}
\end{prop}
\begin{proof}
   Fix a total order $J_1, J_2, \dots, J_k$ of the two-sided $p$-cells of $\W$
   such that $i < j$ implies $J_i \, \overset{p}{\underset{2}{\not\geqslant}} \, J_j$.
   For $0 \leqslant i \leqslant k$ denote by $I_i$ the linear span of $1 \otimes \pkl{x}$
   for $x \in J_s$ with $s \leqslant i$. Then
   \begin{equation}
      \label{eqFilt}
      0 = I_0 \subset I_1 \subset I_2 \subset \dots \subset I_k = \C \W
   \end{equation}
   is a filtration of $\C \W$ by two-sided ideals.
   
   Since $\C \W$ is a semisimple algebra, each simple $\C \W$-module $L$ occurs
   in the left regular representation with multiplicity $\dim(L)$. As \eqref{eqFilt}
   is a filtration by two-sided ideals, there exists an index $1\leqslant i \leqslant k$
   such that $L$ appears with multiplicity $\dim(L)$ in $I_i/I_{i-1}$.
   On the other hand, $I_i/I_{i-1}$ is isomorphic to the direct sum of all
   the cell modules $\weylMod{C}$ for $C$ a left $p$-cell in $J_i$. This implies part (ii)
   of the Proposition. Finally, observe that $L$ occurs exactly once in each 
   of these left cell modules. This proves part (i) and finishes the proof of 
   the Proposition.   
\end{proof}

We extend $\pafp{c}$ to a function $\W \rightarrow \R_{\ge 0}$ mapping $w$ contained
in the left $p$-cell $C$ to $\pafp{c}(C)$. Using the characterization
of the Perron-Frobenius eigenvalue as the largest real eigenvalue, we get
immediately from \cref{thmInv}:

\begin{cor}
   $\pafp{c}$ is constant on two-sided $p$-cells.
\end{cor}

Moreover, we would like to prove:

\begin{conj}
   For $x, y \in W$ we have:
   \[ x \cle{2} y \qquad \Rightarrow \qquad \pafp{c}(x) \geqslant \pafp{c}(y) \]
\end{conj}

Our current techniques do not allow us to prove this as we do not
understand the multiplication of the $p$-canonical basis
well enough. Using Computer calculations we have verified the 
conjecture for finite Weyl groups of rank $\leqslant 4$. 
\Cref{propFamilies} allows us to define:

\begin{defn}
   Let $J \subseteq \W$ be a two-sided $p$-cell and $C \subseteq J$ a left (or right) 
   $p$-cell. We call $L_C$ the $p$-special module of $J$ and the set
   \[\{ L \in \Irr(\W) \; \vert \; [ \weylMod{J}, L ] \neq 0 \}\]
   the $p$-family of $J$.
\end{defn}

In the case of the Kazhdan--Lusztig basis and a left (or right) Kazhdan--Lusztig cell
$C$, the module $L_C$ is a special representation in the sense of Lusztig 
(see \cite[Proposition 9]{KM}). In \cite[Theorem 5.25]{LuChFinRedGrps} Lusztig
shows that for the Kazhdan--Lusztig basis our notion of families coincides with
the original definition. An alternative proof of this result based on the theory 
of primitive ideals can be found in \cite{BVPrimIdealsClassicalGrps, BVPrimIdealsExcepGrps}.

\begin{ex}
   \label{exB2SpecialMods}
   In type $B_2$ we have for the $2$-Kazhdan--Lusztig basis (with the same
   conventions as in \cite{JW}):
   \begin{align*}
      \pkl[2]{sts} &= \kl{sts} + \kl{s} \\
      \pkl[2]{w} &= \kl{w} \text{ for all } w \in \W \setminus \{sts\}
   \end{align*}
   The $2$-special modules in this case are:
   \begin{align*}
      \{\id\} &\rightsquigarrow \text{the sign representation}\\
      \{s\} &\rightsquigarrow \C \text{ where } s  \text{ (resp. } t\text{) acts via } 1 \text{ (resp. } -1\text{)}\\
      \W \setminus \{ \id, s, w_0 \} &\rightsquigarrow \text{the geometric representation}\\
      \{ w_0 \} &\rightsquigarrow \text{the trivial representation}
   \end{align*}
   Thus the irreducible $\C \W$-modules fall into the following $2$-families:
   \[ \{ \sgn \} \cup \{\sgn_s\} \cup \{\sgn_t, \geom\} \cup \{\triv\} \]
   where $\sgn_s$ denotes the special module associated to the two-sided $2$-cell $\{s\}$.
   The difference to Lusztig's families is that the family $\{\sgn_s, \sgn_t, \geom\}$
   splits up into two $2$-families. Moreover, observe that tensoring with the sign
   representation does not give a permutation of the set of $2$-families, whereas
   this is the case for Lusztig's families.
\end{ex}

\section[Left \texorpdfstring{$p$}{p}-Cells Within the Same Two-Sided 
\texorpdfstring{$p$}{p}-Cell]{Study of Left \texorpdfstring{$p$}{p}-Cells 
Within the Same Two-Sided \texorpdfstring{$p$}{p}-Cell}
\label{secPropertyA}

In this section, we will prove that left (or right) $p$-cells within the same
two-sided $p$-cell are incomparable. Again we assume that $W$ is a finite Weyl
group throughout the section. In order to do so, we will need the definition of 
an idempotent two-sided cell which reads in our setting as follows:

\begin{defn}
   Let $J$ be a two-sided $p$-cell. $J$ is called idempotent if there exist
   elements $x, y, z \in J$ such that $\pkl{x}$ occurs with non-trivial
   coefficient in $\pkl{y}\pkl{z}$.
\end{defn}

Even though it might not be obvious at first, we get the following result
(see \cite[Proposition 13 (i)]{KM}):

\begin{lem}
   Each two-sided $p$-cell in $W$ is idempotent.
\end{lem}
\begin{proof}
   Let $J$ be a two-sided $p$-cell. Suppose $J$ is not idempotent.
   The set \[\{\pkl{x} \; \vert \; x \clt{2} J \}\] induces a two-sided ideal 
   in $\heck$ and thus in the semisimple algebra $\C \W$ which we will denote by $I_J$.
   We consider the ideal $I$ spanned by $1 \otimes \pkl{x}$ for $x \in J$ 
   in the semisimple quotient $\C\W / I_J$. As $J$ is not idempotent, 
   $I$ is a non-zero nilpotent ideal, contradicting the semisimplicity
   of the quotient $\C\W / I_J$. (Recall that the Jacobson radical contains all
   nilpotent ideals.)
\end{proof}

Given a left $p$-cell $C$, Kildetoft and Mazorchuk study the set $\mathcal{X}_C$
of all two-sided $p$-cells $J$ such that there exists $x \in J$ satisfying
$\pkl{x} \cdot \heckMod{C} \neq 0$ (or equivalently $\pkl{x} \cdot \weylMod{C} \neq 0$). 
Applying \cite[Propositions 14 and 17]{KM} we get the following properties:
\begin{prop}
   \label{propApex}
   \begin{enumerate}
      \item The set $\mathcal{X}_C$ contains a minimum element, called \emph{the apex 
            of $C$} and denoted $\mathcal{J}(C)$.
      \item For all $x \cge{2} \mathcal{J}(C)$ we have $\pkl{x} \cdot \heckMod{C} \ne 0$.
      \item For any two-sided $p$-cell $I$ and left $p$-cell $C' \subset I$, we have
            $\mathcal{J}(C') = I$.
   \end{enumerate}
\end{prop}

The following result is the main result of this section. In the setting 
of the Kazhdan--Lusztig basis, this result was originally deduced from 
properties of primitive ideals in enveloping algebras (see 
\cite[\S 4]{LuBensonCurtis}). It is a weak form of \cite[P9]{LuUneq}.

\begin{thm}
   If $x \cle{L} y$ and $x \ceq{2} y$, then $x \ceq{L} y$.
\end{thm}

We will break the proof of the theorem into several lemmata,
following along the lines of \cite[Proposition 18 and Corollary 19]{KM}.
Let $J$ be a two-sided $p$-cell and $C \subseteq J$ a left $p$-cell.

\begin{lem}
   \label{lemFPelt}
   The following statements hold:
   \begin{enumerate}
      \item Let $I_J$ be the $\Z[v, v^{-1}]$-span in $\heck$ of all $\pkl{w}$ such that 
            $w \, \overset{p}{\underset{2}{\not\geqslant}} \, J$.
            Then the cell module $\heckMod{C}$ is naturally a $\heck / I_J$-module.
            Similarly, $\weylMod{C}$ is naturally a 
            $\C\W / \C \otimes_{\Z[v, v^{-1}]} I_J$-module.
            
      \item Define $a_{\mathbf{c}} = \sum_{w \in J} c_w \otimes \pkl{w} \in \C \W$.
            Let $M_{\mathbf{c}, C}$ be the matrix giving the action of
            $a_{\mathbf{c}}$ on $\weylMod{C}$. 
            Then $M_{\mathbf{c}, C}$ has positive coefficients.
   \end{enumerate}
\end{lem}
\begin{proof}
   Observe that $I_J$ is a two-sided ideal in $\heck$. Then the first part 
   follows immediately from the definition of the apex of $C$
   and the fact that the apex of $C$ is $J$ (see \cref{propApex} (iii)). 
   It remains to prove the second part.
   
   First, we claim that all columns of $M_{\mathbf{c}, C}$ are non-zero. 
   Suppose that $M_{\mathbf{c}, C}$ has a zero column indexed by $y \in C$. 
   Since the structure coefficients of the $p$-canonical basis are 
   Laurent polynomials with non-negative coefficients (see 
   \cref{propPCanProps} (iii)), this implies 
   \[\pkl{x} \pkl{y} \in \heckMod{\clt{L} C}\] for all $x \in J$. Denote by $I$
   the $\Z[v, v^{-1}]$-span of all $\pkl{x}$ in $\heck / I_J$ for $x \in J$.
   It follows that \[I \cdot \pkl{y} \subset \heckMod{\clt{L} C}\text{.}\] 
   Observe that $I$ is a two-sided ideal in $\heck / I_J$, so that we get:
   \[ I \cdot (\heck / I_J) \cdot \pkl{y} \subset \heckMod{\clt{L} C} \]
   From the transitivity of left cells, it follows  that $\pkl{y}$ generates
   the whole cell module $\heckMod{C}$ under the action of $\heck / I_J$.
   Therefore we must have
   \[ I \cdot \heckMod{C} \subset \heckMod{\clt{L} C} \]
   contradicting the fact that the apex of $C$ is $J$ (see \cref{propApex} (iii)).
   
   Next, we claim that all entries in every column of $M_{\mathbf{c}, C}$ are non-zero.
   Consider the column corresponding to $y \in C$. Let $X \subseteq C$ be the
   set of all $x \in C$ such that the basis element $\pkl{x}$ appears with a 
   non-zero coefficient in $a_{\mathbf{c}} \pkl{y}$ in $\heckMod{C}$. Above we have 
   shown that $X$ is non-empty. Observe that $X$ contains all $x \in C$ such that 
   the basis element $\pkl{x}$ appears with a non-zero coefficient in $I \cdot \pkl{y}$
   in $\heckMod{C}$. Since $I$ is a two-sided ideal in $\heck /I_J$, it thus follows
   that the $\Z[v, v^{-1}]$-span of $\pkl{x}$ for $x \in X$ in $\heckMod{C}$ is 
   invariant under the action of $\heck / I_J$. By the transitivity of the cell 
   module, we see that $X = C$. Therefore, all entries in $M_{\mathbf{c}, C}$ 
   are positive, which concludes the proof.   
\end{proof}

We will continue with the notation of the previous lemma. From \cref{lemFPelt} (ii)
it follows that we can apply the Perron-Frobenius theorem to $M_{\mathbf{c}, C}$.
Let $\lambda$ denote the Perron-Frobenius eigenvalue of $M_{\mathbf{c}, C}$.
From \cref{thmPF} (vi) we deduce that the matrix
\[ M_C \defeq \lim_{m \rightarrow \infty} \frac{M_{\mathbf{c}, C}^m}{\lambda^m} \]
is positive and satisfies $M_C^2 = M_C$. Observe that $M_C$ is the projector
to the one-dimensional $\lambda$-eigenspace of $M_{\mathbf{c}, C}$ and thus 
called \emph{the Perron-Frobenius projector}.
\begin{thm}
   \label{thmIdem}
   We can define an element $e_J \in \C \W/\C \otimes_{\Z[v, v^{-1}]} I_J$ with the following properties:
   \begin{enumerate}
      \item $e_J^2 = e_J$.
      \item $e_J$ can be written as a linear combination of $1 \otimes \pkl{x}$
            for $x \in J$ with positive real coefficients.
      \item $e_J$ acts on $\weylMod{C}$ via $M_C$.
      \item Let $C_1, C_2, \dots, C_k$ be the left $p$-cells in $J$. In any total
            order of the $p$-canonical basis preserving elements of the same left 
            $p$-cell as adjacent elements, we have that $e_J$ acts on $\weylMod{J}$ via 
            \[\begin{pmatrix}
                  M_{C_1} & 0 & \dots & 0 \\
                  0 & M_{C_2} & \ddots & \vdots \\
                  \vdots & \ddots & \ddots & 0 \\
                  0 & \dots & 0 & M_{C_k}
               \end{pmatrix} \text{.}\]
   \end{enumerate}
\end{thm}
\begin{proof}
   First, consider for $m \geqslant 1$ the element 
   \[ \frac{a_{\mathbf{c}}^m}{\lambda^m} = \sum_{x \in J} d_{x, m} \pkl{x} 
      \quad \in \quad \C \W / \C \otimes_{\Z[v, v^{-1}]} I_J \]
   where $d_{x, m} \in \R_{\geqslant 0}$. Since for $x \in J$ the matrix giving 
   the action of $\pkl{x}$ on $\weylMod{C}$ 
   is non-zero by \cref{propApex} (ii) and (iii) and has non-negative coefficients, the 
   sequence $(d_{x, m})_{m \geqslant 1}$ is bounded for all $x \in J$ and 
   thus contains a convergent subsequence.
   
   We have the following identity:
   \begin{align*}
      \frac{a_{\mathbf{c}}^{m+1}}{\lambda^{m+1}} &=  \sum_{x \in J} d_{x, m+1} \pkl{x} \\
         &= \frac{a_{\mathbf{c}}}{\lambda} \cdot \frac{a_{\mathbf{c}}^m}{\lambda^m}
         = \frac{a_{\mathbf{c}}}{\lambda} \cdot \left(\sum_{y \in J} d_{y, m} \pkl{y}\right) \\
         &= \sum_{x, y \in J} \frac{N_{x, y}}{\lambda} d_{y, m} \pkl{x} \\
   \end{align*}
   where $N$ is the matrix giving the action of $a_{\mathbf{c}}$ on the 
   two-sided $p$-cell module $\weylMod{J}$. 
   This implies for all $x \in J$:
   \begin{equation}
      \label{eqCoeffRel}
      d_{x, m+1} = \sum_{y \in J} \frac{N_{x, y}}{\lambda} d_{y, m} 
      = \dots = \sum_{y \in J} \frac{(N^m)_{x, y}}{\lambda^m} d_{y, 1}
   \end{equation}
   Let $C_1, C_2, \dots, C_k$ be an ordering of the left $p$-cells in $J$ such
   that $C_i \cle{L} C_j$ implies $i \leqslant j$. We choose a total order of
   the $p$-canonical basis in $J$ refining the total order of the
   left $p$-cells in $J$. Then observe that $N$ has the following block 
   upper-triangular form 
   \[ N = 
   \begin{pmatrix}
      M_{\mathbf{c}, C_1} & \ast & \dots & \ast \\
      0 & M_{\mathbf{c}, C_2} & \ddots & \vdots \\
      \vdots & \ddots & \ddots & \ast \\
      0 & \dots & 0 & M_{\mathbf{c}, C_k}
   \end{pmatrix}
   \]
   where $M_{\mathbf{c}, C_l}$ is the matrix giving the action of $a_{\mathbf{c}}$
   on the left $p$-cell module $\weylMod{C_l}$ 
   for $1 \leqslant l \leqslant k$.
   Recall that the spectrum of $N$ is the multiset union of the spectra of the
   diagonal block matrices $M_{\mathbf{c}, C_l}$ for $1 \leqslant l \leqslant k$.
   Therefore, \cref{thmPF} implies that $\lambda$ is an eigenvalue with multiplicity
   $k$ for $N$ and that any other eigenvalue $\mu \neq \lambda$ satisfies 
   $\lvert \mu \rvert < \lambda$.
   
   Above we have shown that the sequence $(d_{x, m})_{m \geqslant 1}$ has a convergent
   subsequence for all $x \in J$. This implies that $\lim_{m \rightarrow \infty} 
   \frac{N^m}{\lambda^m}$ exists and consequently that for all $x \in J$ the whole sequence 
   $(d_{x, m})_{m \geqslant 1}$ converges to some $d_x \in \R_{\geqslant 0}$. 
   In addition, it follows that the geometric and algebraic multiplicity of 
   $\lambda$ for $N$ coincide (see \cite[(7.10.33)]{MeMatAna}). 
   
   Next, we define
   \[ e_J \defeq \sum_{x \in J} d_x \pkl{x} \text{.} \]
   It follows immediately that $e_J^2 = e_J$ as $e_J$ is the projector to the 
   $\lambda$-eigenspace of $N$ and that $e_J$ acts via $M_C$ on $\weylMod{C}$.
   This concludes the proof of (i) and (iii).
   
   We will prove part (iv) next. Observe that the action of $e_J$ on the
   two-sided $p$-cell module $\weylMod{J}$ in the $p$-canonical basis is given by
   the following matrix:
   \[ N_J \defeq \lim_{m \rightarrow \infty} \frac{N^m}{\lambda^m} = 
   \begin{pmatrix}
      M_{C_1} & \ast & \dots & \ast \\
      0 & M_{C_2} & \ddots & \vdots \\
      \vdots & \ddots & \ddots & \ast \\
      0 & \dots & 0 & M_{C_k}
   \end{pmatrix}
   \]
   which is block upper-triangular with positive, idempotent matrices on the 
   diagonal (by the considerations preceding the proposition). Since $N_J$ is
   a non-negative idempotent matrix, we may apply \cite[Theorem 2]{FlGrpsNonNegMats}
   to get that all off-diagonal blocks are $0$. (Actually, the easy special case
   where neither a row nor a column of the non-negative idempotent matrix is zero 
   suffices.) Therefore, $N_J$ is the direct sum of the positive, idempotent
   matrices $M_{C_i}$ for $1 \leqslant i \leqslant k$. This completes
   the proof of part (iv).

   It remains to prove part (ii). From \eqref{eqCoeffRel}, we deduce that:
   \[ d_x = \sum_{y \in J} N_{x, y} c_{y, 1} = \sum_{1\leqslant l \leqslant k} 
      \sum_{y \in C_l} (M_{C_l})_{x, y} \frac{c_y}{\lambda} > 0 \]
   This concludes the proof. Observe that $e_J$ projects in each cell 
   module $\weylMod{C_l}$ for $1 \leqslant l \leqslant k$ to the the unique (up to scalar) 
   non-negative eigenvector of $M_{C_l}$ by \cref{thmPF} (iv) and (v).
   Moreover, the image of $e_J$ in $\weylMod{J}$ gives a positive
   eigenvector for $N_J$ for the eigenvalue $1$ which is not unique
   though as $1$ is a semisimple eigenvalue of multiplicity $k$.
\end{proof}

Recall that the coefficients ${}^p \mu_{x, y}^z$ for $x, y, z \in W$ are
the structure coefficients of the $p$-canonical basis (see \cref{propPCanProps} (iii)). 
The reader should compare the following corollary with \cite[P8]{LuUneq} and
keep in mind that $\gamma_{x,y,z}$ is the coefficient in front of the highest
(or lowest) power of $v$ in ${}^0 \mu_{x, y}^{z^{-1}}$. One would obtain the
third cell equivalence if cyclicity (see \cite[P7]{LuUneq}) held.
\begin{cor}
   Let $x, y, z \in J$. If ${}^p \mu_{x, y}^z$ is non-zero, then
   $y \ceq{L} z$ and $x \ceq{R} z$. In other words, $z$ lies in the
   intersection of the right $p$-cell of $x$ with the left $p$-cell of $y$.
   Moreover, any $\pkl{z}$ for $z \in J$ occurs with non-trivial coefficient in a
   product $\pkl{x} \cdot \pkl{y}$ with $x, y \in J$.
\end{cor}
\begin{proof}
   \cref{thmIdem} (ii) shows that $e_J$ is a positive linear combination
   of all $p$-canonical basis elements indexed by elements in $J$.
   The block-diagonal form of the matrix $N_J$ giving the action of $e_J$ on 
   the two-sided $p$-cell module $\weylMod{J}$ (see \cref{thmIdem} (iv)) 
   demonstrates that ${}^p \mu_{x, y}^z$ non-zero implies $y \ceq{L} z$.
   
   To see that any $\pkl{z}$ for $z \in J$ occurs with non-trivial coefficient 
   in a product $\pkl{x} \cdot \pkl{y}$ with $x, y \in J$, recall
   that all the diagonal blocs in $N_J$ have positive coefficients. 
   
   Recall that $J^{-1}$ is a two-sided $p$-cell as well. Applying the 
   $\Z[v, v^{-1}]$-linear anti-involution $\iota$, we obtain 
   ${}^p \mu_{x, y}^z = {}^p \mu_{y^{-1}, x^{-1}}^{z^{-1}}$. Combined with the first 
   part of the corollary applied to $J^{-1}$ we see that ${}^p \mu_{x, y}^z$ non-zero
   gives $x^{-1} \ceq{L} z^{-1}$ which is equivalent to $x \ceq{R} z$.
\end{proof}

\begin{cor}
   \label{corIneqLCells}
   The left $p$-cells within $J$ are incomparable with respect to
   the left $p$-cell preorder. In other words, we have for $x, y \in J$:
   \[ x \cle{L} y \quad \Rightarrow \quad x \ceq{L} y \]
\end{cor}
\begin{proof}
   Let $I$ be the $\Z[v, v^{-1}]$-span of all $\pkl{z}$ in $\heck / I_J$ for $z \in J$.
   If $x \cle{L} y$, then there exists $h \in \heck$ such that $\pkl{x}$ occurs
   with non-trivial coefficient in $h \cdot \pkl{y}$. From the previous corollary it follows
   that there are $v, w \in J$ such that ${}^p \mu_{v, w}^y$ is non-zero.
   Thus, $\pkl{x}$ occurs with non-trivial coefficient in $h \pkl{v} \pkl{w}$.
   Since $I$ is a two-sided ideal in $\heck / I_J$, it follows that the image of
   $h \pkl{v} \pkl{w}$ in $\heck / I_J$ lies in $I$. Rewriting 
   $(h \pkl{v}) \pkl{w}$ in $\heck / I_J$ in terms of the $p$-canonical basis 
   and applying the previous corollary shows that $x$ lies
   in the same left $p$-cell as $w$ which in turn lies in the same left cell as $y$.
\end{proof}

\begin{cor}
   Two-sided $p$-cells are the smallest subsets that are at the same time
   a union of left $p$-cells and one of right $p$-cells.
\end{cor}

\begin{ex}
   Let us illustrate the results of this section in type $B_2$ continuing 
   \cref{exB2SpecialMods}. Let us consider the most interesting two-sided
   $2$-cell $J = \{ 2, 12, 212, 21, 121 \}$ which decomposes into two
   left $2$-cells
   \[ J = \{2, 12, 212 \} \cup \{21, 121\} \]
   and the element $a_{\mathbf{c}} = \sum_{x \in J} \pkl{x}$ (i.e. all
   constants are chosen to be $1$). The matrix giving the action of $a_{\mathbf{c}}$
   on $\weylMod{J}$ looks as follows
   \[ N \defeq \begin{pmatrix}
         3 & 4 & 3 &  &  \\
         4 & 6 & 4 &  &  \\
         3 & 4 & 3 &  &  \\
          &  &  & 6 & 4 \\
          &  &  & 8 & 6
      \end{pmatrix} \] 
   where we omitted the $0$ entries. $N$ has as Perron-Frobenius eigenvalue 
   $\lambda = 6 + 4\sqrt{2}$. The eigenvectors to this eigenvalue for 
   multiplication with $N$ on the right are given by
   \[ \hat{v}_1^{\tr} \defeq \begin{pmatrix}
         1 & \sqrt{2} & 1 & 0 & 0
      \end{pmatrix}
      \text{ and }
      \hat{v}_2^{\tr} \defeq \begin{pmatrix}
         0 & 0 & 0 & \sqrt{2} & 1
      \end{pmatrix} \text{.} \]
   Similarly, the eigenvectors to this eigenvalue for multiplication with $N$ on the left 
   are given by $v_1 \defeq \hat{v}_1$ and 
   \[ v_2 \defeq
      \begin{pmatrix}
         0 \\ 0 \\ 0 \\ \sqrt{2} \\ 2
      \end{pmatrix} \text{.} \]
   After normalization, these eigenvectors allow us to easily calculate $N_J$:
   \[ N_J \defeq \lim_{m \rightarrow \infty} \frac{N^m}{\lambda^m} = 
   \begin{pmatrix}
      \frac{1}{4} & \frac{\sqrt{2}}{4} & \frac{1}{4} &  &  \\[3pt]
      \frac{\sqrt{2}}{4} & \frac{1}{2} & \frac{\sqrt{2}}{4} &  &  \\[3pt]
      \frac{1}{4} & \frac{\sqrt{2}}{4} & \frac{1}{4} &  &  \\[3pt]
       &  &  & \frac{1}{2} & \frac{\sqrt{2}}{4} \\[3pt]
       &  &  & \frac{\sqrt{2}}{2} & \frac{1}{2}
   \end{pmatrix} = \frac{1}{4} (v_1 \hat{v}_1^{\tr} + v_2 \hat{v}_2^{\tr}) \]
   For the idempotent $e_J$ we get:
   \begin{align*}
      e_J \defeq &\frac{2 - \sqrt{2}}{8} \otimes (\pkl{2} + \pkl{212} + \pkl{121}) \\
           + &\frac{\sqrt{2} - 1}{4} \otimes (\pkl{21} + \pkl{12})
           \quad \in \quad  \C\W / \C \otimes_{\Z[v, v^{-1}]} \heckMod{\cle{2} 1212}
   \end{align*}
\end{ex}

\section{\texorpdfstring{$p$}{p}-Cells for Symmetric Groups}
\label{secPKLCellular}

\subsection{Two-Sided \texorpdfstring{$p$}{p}-Cell Preorder and the Dominance Order}

The Robinson-Schensted correspondence (see \cite[\S A.3.3]{BBCoxGrps} or 
\cite[\S4.1]{FuTableaux}) gives a bijection between the symmetric group $S_n$ and pairs 
of standard tableaux of the same shape with $n$ boxes. The row-bumping algorithm gives 
a way to explicitly calculate the image $(P(w), Q(w))$ of $w \in S_n$ under the 
Robinson-Schensted correspondence.

Throughout this section we assume that we used the Cartan matrix in finite type 
$A_{n-1}$ as input. In this case the Coxeter system $(W, S)$ can be identified with 
the pair $(S_n, \{s_1, \dots, s_{n-1}\})$ consisting of the symmetric group 
together with the set of simple transpositions.

$p$-Cells for symmetric groups admit a beautiful description in terms of the
Robinson-Schensted correspondence (see \cite[Theorem 4.33]{JeABC}):
\begin{thm}
   \label{thmPCellsTypeA}
   For $x, y \in S_n$ we have:
   \begin{align*}
      x \ceq{L} y &\Leftrightarrow Q(x) = Q(y) \\
      x \ceq{R} y &\Leftrightarrow P(x) = P(y) \\
      x \ceq{LR} y &\Leftrightarrow Q(x) \text{ and } Q(y) \text{ have the same shape}
   \end{align*}
   In particular, Kazhdan--Lusztig cells and $p$-cells of $S_n$ coincide.
\end{thm}

The last result shows that partitions of $n$ correspond to two-sided cells of 
$S_n$. Recall the definition of the dominance order on the set of partitions of $n$:
\begin{defn}
   Let $\lambda = (\lambda_1, \lambda_2, \dots), \mu=(\mu_1, \mu_2, \dots)$ be 
   two partitions of $n$. We have $\lambda <= \mu$ in the dominance order
   if and only if
   \[ \sum_{i=1}^k \lambda_i \leqslant \sum_{i=1}^k \mu_i \]
   for all $k \geqslant 1$.
\end{defn}

The following fundamental property of the dominance order (see 
\cite[Proposition 2.8]{BrPartitions}) will be useful for us:
\begin{prop}
   \label{propTr}
   For two partitions $\lambda$, $\mu$ of $n$ the following holds:
   \[ \lambda \leqslant \mu \quad \Leftrightarrow \quad \lambda^T \geqslant \mu^T \]
   where $\lambda^T$ denotes the \emph{conjugate or transpose partition of $\lambda$}.
\end{prop}

We will not explicitly need the following definition, but we give it for the sake
of completeness (see \cite[\S A3.7 and \S A3.8]{BBCoxGrps} for details):
\begin{defn}
   Given a standard tableau $T$, the combinatorial algorithm called
   \emph{evacuation} proceeds as follows:
   
   First, delete the entry $1$ and perform a backward slide on the cell that
   contained it. Then repeat this for the entry $2$, etc.
   Finally, the tableau that records in reverse the order in which the cells 
   of $T$ have been vacated is called the \emph{evacuation of $T$} and denoted
   by $e(T)$.
\end{defn}

For our next result, we will use the following compatibility of the Robinson-Schensted
correspondence with the multiplication of the longest element $w_0 \in S_n$
(see \cite[Fact A3.9.1]{BBCoxGrps} or \cite[Theorem D]{KnArtIII} for a proof):
\begin{thm}
   \label{thmLongestElt}
   If $x \in S_n$ corresponds to $(P, Q)$ under the Robinson-Schensted correspondence,
   then we have under the Robinson-Schensted correspondence that
   \begin{enumerate}
      \item $x w_0$ corresponds to $(P^T, e(Q)^T)$,
      \item $w_0 x$ corresponds to $(e(P)^T, Q^T)$ and
      \item $w_0 x w_0$ corresponds to $(e(P), e(Q))$
   \end{enumerate}
   where $e(P)$ is the \emph{evacuation} of $P$.
\end{thm}

To simplify notation, we will denote by $J_{\lambda} \subseteq S_n$  the two-sided
cell corresponding to a partition $\lambda$ of $n$. In order to prove the main result, 
we will rely on the following observation which is motivated by 
\cite[Proposition A.2.1]{W1}:
\begin{prop}
   \label{propWeakOrderGenDomOrder}
   The dominance order on partitions of $n$ is generated by the weak Bruhat order
   under the Robinson-Schensted correspondence.
\end{prop}
\begin{proof}
   Let $\lambda, \mu$ be two partitions of $n$ such that $\lambda \geqslant \mu$.
   We will prove the statement by induction on the length of a maximal chain between 
   $\lambda$ and $\mu$. For $\lambda = \mu$ there is nothing to show. 
   
   For the induction step, suppose $\lambda > \mu$. It is enough to show that there 
   exists $s \in S$, $x \in J_{\mu}$ and a partition $\nu$ of $n$ such that
   \begin{itemize}
      \item $xs < x$,
      \item $xs \in J_{\nu}$ and
      \item $\lambda \geqslant \nu > \mu$.
   \end{itemize}
   
   The idea is to obtain $\nu$ from $\mu$ by applying a single raising
   operation. It is a well-known fact that the raising operations generate the
   dominance order (see for example \cite[(1.16)]{MacSymmFuncs}).
   
   Let $i \in \N$ be minimal such that $\lambda_i > \mu_i$ and observe that
   it is the first index for which the parts of $\lambda$ and $\mu$ differ.
   It follows $\mu_{i-1} = \lambda_{i-1} \geqslant \lambda_i > \mu_i$
   and $\mu_{i+1} \ne 0$ as $\mu$ and $\lambda$ are both partitions of $n$.
   
   We will distinguish two cases:
   \begin{description}
   \item[1. Case:] $\mu_{i+1} > \mu_{i+2}$.
   Define the partition $\nu$ as follows:
   \[ v_j \defeq \begin{cases}
      \mu_i + 1 & \text{ if } j = i\text{,} \\
      \mu_{i+1} - 1 & \text { if } j = i+1\text{,} \\
      \mu_j & \text{ otherwise.}
   \end{cases} \]
   It follows immediately that $\nu$ is a partition of $n$ satisfying 
   $\lambda \geqslant \nu > \mu$.
   
   Our convention of drawing partitions is that $\nu_i$ gives the number
   of boxes in row $i$. With this convention in mind, divide $\nu$ into
   three parts:
   \begin{itemize}
   \item the first $i - 1$ rows,
   \item the $i$-th and $(i+1)$-st row and
   \item the remaining rows.
   \end{itemize}
   
   Let $T$ be the standard tableau of shape $\nu$ that is column superstandard 
   in each of the three pieces (of course up to shift so that the numbers do not 
   repeat). Suppose that there are $k$ boxes in the first $i-1$ rows. Then
   rows $i$ and $i+1$ of $T$ look as follows:
   
   \begin{center}
      \resizebox{\linewidth}{!}{%
      \begin{tikzpicture}[scale=1.5]
         \draw (0,0) -- (8,0);
         \draw (0,-1) -- (8,-1);
         \draw (0, -2) -- (4, -2);
         
         \draw (0, -2.1) -- (0, 0.1);
         \draw (1, -2.1) -- (1, 0.1);
         \draw (2, -2.1) -- (2, 0.1);
         \draw (3, -2.1) -- (3, 0.1);
         \draw (4, -2.1) -- (4, 0.1);
         
         \draw (5, -1) -- (5, 0.1);
         \draw (6, -1) -- (6, 0.1);
         \draw (7, -1) -- (7, 0.1);
         \draw (8, -1) -- (8, 0.1);
         
         \node[labelling] (vkp2) at (0.5, -1.5) {\small$k+2$};
         \node[labelling] (vkp4) at (1.5, -1.5) {\small$k+4$};
         \node[labelling] (vkdots3) at (2.5, -1.5) {\small$\cdots$};
         \node[labelling] (vkpnu) at (3.5, -1.5) {\small$k + 2\nu_{i+1}$};
         
         \node[labelling] (vkp1) at (0.5, -0.5) {\small$k + 1$};
         \node[labelling] (vkp3) at (1.5, -0.5) {\small$k + 3$};
         \node[labelling] (vkdots) at (2.5, -0.5) {\small$\cdots$};
         \node[labelling, align=center] (vkpnum1) at (3.5, -0.5) {\small$k + 2\nu_{i+1}$\\ \small$-1$};
         \node[labelling, align=center] (vkpnup1) at (4.5, -0.5) {\small$k+2\nu_{i+1}$\\ \small$+1$};
         \node[labelling, align=center] (vkpnup2) at (5.5, -0.5) {\small$k + 2\nu_{i+1}$\\ \small $+2$};
         \node[labelling] (vkdots2) at (6.5, -0.5) {\small$\cdots$};
         \node[labelling, align=center] (vkpnupnu) at (7.5, -0.5) {\small$k + \nu_{i+1}$\\ \small$+\nu_{i}$};
      \end{tikzpicture}}
   \end{center}
   
   Define a permutation $y \in S_n$ as follows: Let $(y(1), y(2), \dots, y(n))$ be 
   the reading word of $T$ (i.e. the sequence obtained by reading the entries
   of $T$ from bottom to top, left to right). Then the sequence contains the following
   piece for part $2$:
   \begin{align*}
      (\dots,\,&k+2,\,k+4,\,\dots,\,k+2\nu_{i+1},\,k+1,\,k+3,\,\dots,\,k+2\nu_{i+1}-1, \\
      &k+2\nu_{i+1}+1,\,k+2\nu_{i+1}+2,\,\dots,\,k+\nu_{i+1}+\nu_i,\,\dots)
   \end{align*}
   Observe that $k+\nu_{i+1}+\nu_i$ and $k+\nu_{i+1}+\nu_i-1$ are in neighbouring
   positions and let $s \in S$ be the simple transposition swapping these two positions.
   
   It follows $ys > y$ as $s$ introduces a new inversion. We have that $ys$ lies
   in $J_{\mu}$ because when applying the row-bumping algorithm the entry 
   $k+\nu_i+\nu_{i+1}$ is bumped one row down by $k+\nu_i + \nu_{i+1}-1$
   and thus ends up in row $i+1$. This finishes the proof in the first case
   by setting $x \defeq ys$.
    
   \item[2. Case:] $\mu_{i+1} = \mu_{i+2}$. Let $m \in \N$ be maximal such that
   $\mu_m = \mu_{i+1}$. Define the partition $\nu$ as follows:
   \[ v_j \defeq \begin{cases}
      \mu_i + 1 & \text{ if } j = i\text{ and }\mu_i = \mu_{i+1}\text{,} \\
      \mu_{i+1} + 1 & \text { if } j = i+1\text{ and }\mu_i > \mu_{i+1}\text{,} \\
      \mu_m - 1 & \text{ if } j = m\text{,}\\
      \mu_j & \text{ otherwise.}
   \end{cases} \]
   It follows that $\nu$ is a partition of $n$ satisfying $\lambda \geqslant \nu > \mu$. 
   Indeed, observe that the smallest index $k \in \N$ such that 
   $\sum_{l = 1}^k \lambda_l = \sum_{l = 1}^k \mu_l$ satisfies $k \geqslant m$
   in order to see $\lambda \geqslant \nu$. 
   
   Next, observe that $\nu$ and $\mu$ differ in exactly two adjacent columns $l$ and $l+1$.
   By \cref{propTr}, we have $\nu^T < \mu^T$ and $\mu^T$ is obtained from $\nu^T$
   by raising a box from row $l+1$ to row $l$. We will define the permutation
   $y \in S_n$ and the simple transposition $s \in S$ as before, with the
   only difference that $T$ will be of shape $\mu^T$ instead of $\nu$. 
   It follows that $ys$ lies in $J_{\nu^T}$ and satisfies $ys > y$.
   Define a simple reflection $t \defeq w_0 s w_0$. \Cref{thmLongestElt} shows 
   that $x \defeq yw_0$ lies in $J_{\mu}$ and $x t = y s w_0$ lies in $J_{\nu}$.
   Moreover, we have $ x t = y s w_0 < y w_0 = x$. This finishes
   the proof of the proposition.\qedhere
   \end{description}
\end{proof}

The goal of this section is to prove the following result:
\begin{thm}
   \label{thmOrders}
   Let $J, J' \subseteq S_n$ be two-sided cells that correspond to the partitions
   $\lambda, \lambda'$ of $n$ respectively. Then we have the following:
   \[ J \cle{2} J' \Leftrightarrow \lambda \leqslant \lambda' \]
   In particular, the two-sided $p$-cell preorder coincides with the Kazhdan--Lusztig
   two-sided cell preorder for all primes $p$.
\end{thm}
\begin{proof}
   \begin{description}[leftmargin=0pt]
   \item[$\Rightarrow$] In this case, we have $h, h' \in \heck$ such that
                       for some $w \in J$ the $p$-canonical basis element
                       $\pkl{w}$ occurs with non-trivial coefficient in $h \pkl{w_{J'}} h'$
                       where $w_{J'}$ is the longest element in a standard parabolic
                       subgroup contained in $J'$. Since the corresponding Schubert
                       variety is smooth, we have $\pkl{w_{J'}} = \kl{w_{J'}}$.
                       This implies
                       \[ w \cle[0]{2} w_{J'} \text{ and thus } J \cle[0]{2} J'\text{.} \]
                       Applying the characterization of the two-sided Kazhdan--Lusztig
                       preorder in terms of the dominance order (see 
                       \cite[Theorem 5.1]{GeMurphyBasis}), we get 
                       $\lambda \leqslant \lambda'$. 
                     
                       It is important to note that Geck uses the opposite dominance 
                       order because his connection to two-sided Kazhdan--Lusztig cells 
                       is not based on the Robinson-Schensted correspondence, but on
                       leading matrix coefficients of irreducible representations of 
                       $\heck \otimes_{\Z[v, v^{-1}]} \R(v)$. It depends on the
                       parametrization of irreducible representations of $S_n$
                       in a way that is compatible with induced sign characters from
                       Young subgroups (see \cite[Example 3.8]{GeMurphyBasis}).
                       In his correspondence, he introduces the transpose of a partition
                       which by \cref{propTr} inverses the dominance order
                       (compare \cite[Corollary 5.6]{GeMurphyBasis} with 
                       \cref{thmPCellsTypeA}).
                       
   \item[$\Leftarrow$] By induction on the maximal chain between $\lambda$ and $\mu$
                       it is enough to prove $J \clt{2} J'$ in the case where 
                       $\lambda < \lambda'$ are adjacent in the dominance order (so that 
                       $\lambda \leqslant \mu \leqslant \lambda'$ implies $\mu = \lambda$
                       or $\mu = \lambda'$). In the proof of 
                       \cref{propWeakOrderGenDomOrder} we have shown that in this case
                       there exists $x \in J$, $s \in S$ such that $xs < x$ and $xs \in J'$.
                       Due to $(xs)s = x > xs$ it follows that $\pkl{x}$ occurs with 
                       non-trivial coefficient in $\pkl{xs} \kl{s}$. This implies 
                       $x \clt{R} xs$ and thus $J \clt{2} J'$ finishing the proof.\qedhere
   \end{description}
\end{proof}

The last result in the setting of the Kazhdan--Lusztig basis is a special
case of a conjecture made independently by Lusztig in 
\cite[\S 10.8]{LuCellsInAffWeylGrpsIV} and Vogan in \cite{VoOrbitMethod}. Proofs can 
be found in \cite[2.13.1]{DPSCellsQSChur} and \cite[Theorem 5.1]{GeMurphyBasis} 
for symmetric groups, in \cite{ShPartialOrderOnCells} for affine Weyl groups of 
type $\widetilde{A}_n$ or  of rank $\leqslant 4$ and in \cite[Theorem 4 b)]{BLuBij} 
in general.

It remains an open question to relate the left (or right) $p$-cell preorder and 
the Kazhdan--Lusztig left (or right) cell preorder. It should also be noted that
a charaterization of either of these preorders on standard tableau under
the Robinson-Schensted correspondence is not known.

\subsection{Cellularity of the \texorpdfstring{$p$}{p}-Canonical Basis}

First, we will recall the definition of a cellular algebra given in 
\cite[Definition 1.1]{GLCellAlgs}. Let $R$ be a fixed commutative ring with 
unit and let $A$ be an $R$-algebra.

\begin{defn}
   A \emph{cell datum} for $A$ is a quadruple $(\Lambda, \ast, M, C)$ consisting of:
   \begin{itemize}
      \item A finite partially ordered set $\Lambda$,
      \item An $R$-linear anti-involution $\ast$ of $A$,
      \item For every $\lambda \in \Lambda$ a finite, non-empty set $M(\lambda)$ 
            of indices, and
      \item An injective map $C: \bigcup_{\lambda \in \Lambda} M(\lambda) \times 
            M(\lambda) \rightarrow A$. If $\lambda \in \Lambda$ and $S, T \in 
            M(\lambda)$ write $C^{\lambda}_{S, T} = C(S, T) \in A$.
   \end{itemize}
   and satisfying the following conditions:
   \begin{enumerate}
      \item The image of $C$ gives an $R$-basis of $A$.
      \item $(C^{\lambda}_{S, T})^{\ast} = C^{\lambda}_{T, S}$.
      \item If $\lambda \in \Lambda$ and $S, T \in M(\lambda)$ then we have for any element
            $a \in A$
            \[ a C^{\lambda}_{S, T} = \sum_{S' \in M(\lambda)} r_a(S', S) 
               C^{\lambda}_{S', T} \quad \text{(mod } A(< \lambda) ) \]
            where the coefficients $r_a(S', S) \in R$ are independent of $T$ and where
            $A(< \lambda)$ denotes the $R$-submodule of $A$ generated by 
            $\{C^{\mu}_{V, W} \; \vert \; \mu < \lambda, V, W \in M(\mu)\}$.
   \end{enumerate}
\end{defn}

Examples of algebras that can be equipped with a cell datum include
matrix rings $\Mat_{d \times d}(R)$, $R[x]/(x^n)$, the Temperley Lieb 
algebra and many more (see \cite[\S1, \S4-6]{GLCellAlgs} for the details).

For the rest of this section we assume again that we used the Cartan matrix in 
finite type $A_{n-1}$ as input for the Hecke category.

The goal of this section is to extend the $p$-canonical basis to a cell datum of
$\heck$. The cell datum will be chosen as follows:
\begin{itemize}
   \item Let $\Lambda$ be the set of partitions of $n$, equipped with the dominance
         order.
   \item Let $\ast$ be the $\Z[v, v^{-1}]$-linear anti-involution $\iota$ of $\heck$.
   \item For $\lambda \in \Lambda$, let $M(\lambda)$ be the set of standard
         tableaux of shape $\lambda$.
   \item If $w \in S_n$ corresponds to $(P, Q)$ of shape $\lambda$ under the 
         Robinson-Schensted correspondence, then define $C$ to map $(P, Q)$
         to $\pkl{w}$. Since the Robinson-Schensted correspondence gives a bijection
         between $S_n$ and the set of pairs of standard tableaux of the same shape
         with $n$ boxes, the map $C$ is obviously injective.         
\end{itemize}

Before we can prove the main theorem, we will need the following 
classical result about the Robinson-Schensted correspondence (see 
\cite[\S4.1, Corollary to Symmetry Theorem]{FuTableaux}):

\begin{thmlab}[Symmetry Theorem for $S_n$]
   \label{thmSymm}
   If $w \in S_n$ corresponds to $(P(w), Q(w))$, then $w^{-1}$ corresponds to
   $(Q(w), P(w))$ under the Robinson-Schensted correspondence.
\end{thmlab}

Moreover, we need the following classical result by Knuth (see \cite[Theorem 6]{Kn}):
\begin{thm}
   \label{thmKEqPSymbols}
   Let $x, y \in S_n$. Then $x$ and $y$ are Knuth equivalent if and only if
   $P(x) = P(y)$.
\end{thm}

In addition, we will require the following description of descent sets under
the Robinson-Schensted correspondence (see \cite[Fact A3.4.1]{BBCoxGrps} and 
\cite[Proposition 2.7.1]{W1} for a proof):
\begin{defn}
   Let $P$ be a standard tableau. The \emph{tableau descent set of $P$}, denoted
   by $\desc{D}(P)$, is the set of integers $i > 0$ for which $i+1$ lies strictly
   below and weakly to the left of $i$ in $P$.
\end{defn}
\begin{lem}
   \label{lemDescSets}
   For $w \in S_n$ we have:
   \begin{enumerate}
      \item $s_i \in \desc{L}(w) \Leftrightarrow i \in \desc{D}(P(w))$
      \item $s_i \in \desc{R}(w) \Leftrightarrow i \in \desc{D}(Q(w))$
   \end{enumerate}
\end{lem}

As last ingredient, we need to show that when acting on a left cell module the
resulting structure coefficients are independent of the $Q$-symbols involved.
\begin{lem}
   \label{lemStructCoeffs}
   Let $x, x', y, y' \in S_n$ be such that $x \ceq{L} y$, $x' \ceq{L} y'$, 
   $P(x) = P(x')$, and $P(y) = P(y')$. Then we have for all $s \in S$ the following:
   \[ {}^p \mu_{s, x}^y = {}^p \mu_{s, x'}^{y'} \]
   In particular, we may introduce the notation 
   ${}^p r_s(P(y), P(x)) \defeq {}^p \mu_{s, x}^y$ as this coefficient does not depend
   on the $Q$-symbols of $x$ and $y$.
\end{lem}
\begin{proof}
   First, observe that $\desc{L}(x) = \desc{L}(x')$ by \cref{lemDescSets} as the 
   $P$-symbols of $x$ and $x'$ coincide. We will consider the only interesting
   case where $s \notin \desc{L}(x)$.
   
   By \cref{thmKEqPSymbols}, there exists a sequence of elementary Knuth 
   operations relating $x$ and $x'$ as their $P$-symbols coincide.
   Each elementary Knuth operation corresponds to a right star operation with respect
   to a rank $2$ standard parabolic subgroup (see \cite[Lemma 4.30]{JeABC}) and thus we
   will use the corresponding terms interchangeably in the following.
   
   Next, we claim that we may apply the whole sequence of right star operations to $y$
   as well. Since $x$ and $y$ lie in the same left cell, we have 
   $\desc{R}(x) = \desc{R}(y)$ by combining \cref{thmPCellsTypeA} and \cref{lemDescSets}.
   Therefore, the first right star operation can be applied to $y$ as well. 
   By \cite[Theorem 4.13]{JeABC}, applying the first right star operation to $x$ and to 
   $y$ gives two elements that still lie in the same left cell. Therefore, we can 
   repeat the argument to see that we can apply the whole sequence of elementary 
   Knuth transformations to $y$ to obtain an element $y''$.
   
   Since Knuth equivalent permutations have the same $P$-symbol (see 
   \cref{thmKEqPSymbols}), we have $P(y') = P(y) = P(y'')$. Our argument 
   also shows that $x'$ and $y''$ lie in the same left cell. Therefore, 
   we have $Q(y') = Q(x') = Q(y'')$ by \cref{thmPCellsTypeA}. It follows 
   that $y' = y''$ as both their $P$ and $Q$-symbols coincide.
   
   The result now follows from repeated application of \cite[Corollary 4.10]{JeABC}.
\end{proof}

Finally, we can prove our main result of this section:

\begin{thm}
   The quadruple $(\Lambda, \ast, M, C)$ gives a cell datum for $\heck$.
\end{thm}
\begin{proof}
   (i) follows from the fact that the $p$-canonical basis is a basis for the 
   $\Z[v, v^{-1}]$-algebra $\heck$ and that the Robinson-Schensted correspondence
   gives a bijection between the indexing sets involved.
   
   Combining $\iota(\pkl{x}) = \pkl{x^{-1}}$ for all $x \in \W$ 
   (see \cite[Proposition 2.5 (iv)]{JeABC}) with \cref{thmSymm} gives (ii).
   
   It remains to prove condition (iii). First, let us recall that we have
   for $s_i \in S$ and $x \in W$:
   \begin{equation}
      \label{eqMultForm}
      \kl{s_i} \pkl{x} = 
      \begin{cases}
         (v + v^{-1}) \pkl{x} & \text{if } s_i \in \desc{L}(x)\text{,}\\
         \sum_{\substack{y \leqslant s_i x\\s_i y < y}} {}^p \mu_{s_i, x}^y \pkl{y} &
         \text{otherwise.}
      \end{cases}
   \end{equation}
   
   Let $J$ be the two-sided cell of $x$. We want to reduce the multiplication
   formula \eqref{eqMultForm} modulo $\heckMod{\clt{2} J}$. If $\pkl{y}$ for some
   $y \in \W$ occurs with non-zero coefficient in \eqref{eqMultForm}, then
   we have $y \cle{L} x$. \cref{corIneqLCells} implies that if $\pkl{y}$ does not
   lie in $\heckMod{\clt{2} J}$, then we have $y \ceq{L} x$. In particular,
   $x$ and $y$ have the same $Q$-symbol under the Robinson-Schensted correspondence
   by \cref{thmPCellsTypeA}.
   
   Let $\lambda$ be the partition of $n$ that corresponds to the two-sided cell $J$.
   By \cref{thmOrders}, the $\Z[v, v^{-1}]$-submodule $\heckMod{\clt{2} J}$ coincides
   with $\heck(< \lambda)$.
   
   \Cref{lemDescSets} shows that for the left descent set of a permutation coincides with
   the left descent set of its $P$-symbol. Suppose that $x$ corresponds to $(P, Q)$ 
   under the Robinson-Schensted correspondence to simplify notation. 
   
   Combining the arguments above with \cref{lemStructCoeffs} we may rewrite and 
   reduce \eqref{eqMultForm} modulo $\heck(< \lambda)$ as follows:   
   \begin{equation}
      \kl{s_i} C_{P, Q}^{\lambda} = 
      \begin{cases}
         (v + v^{-1}) C_{P, Q}^{\lambda} & \text{if } i \in \desc{D}(P)\text{,}\\
         \sum_{\substack{P' \in M(\lambda)\\i \in \desc{D}(P')}} {}^p r_{s_i}(P', P) 
               C_{P', Q}^{\lambda} & \text{otherwise.}
      \end{cases}
      \quad (\text{mod }H(< \lambda))
   \end{equation}
   This finishes the proof of (iii).
\end{proof}

\printbibliography

\end{document}